\documentclass[leqno,11pt]{article}%
\usepackage{amsfonts}
\usepackage{amsmath}
\usepackage{amssymb}
\usepackage{graphicx}%
\providecommand{\U}[1]{\protect\rule{.1in}{.1in}}
\newtheorem{theorem}{Theorem}

\newtheorem{problem}[theorem]{Problem}
\newtheorem{proposition}[theorem]{Proposition}
\newtheorem{remark}[theorem]{Remark}

\newenvironment{proof}[1][Proof]{\noindent\textbf{#1.} }{\ \rule{0.5em}{0.5em}}

\begin{document}

\title{{\LARGE Numerical Schemes for Multivalued Backward Stochastic Differential
Systems}$\smallskip$}
\author{{\large Lucian Maticiuc}$^{\ast}${\large , Eduard Rotenstein}$^{\ast}$%
\bigskip\bigskip\\\textit{Faculty of Mathematics, \textquotedblleft Alexandru Ioan
Cuza\textquotedblright\ University,}\\\textit{Carol I Blvd., no.~9, 700506, Ia\c{s}i, Rom\^{a}nia.}\medskip}
\date{}
\maketitle

\begin{abstract}
We define some approximation schemes for different kinds of generalized
backward stochastic differential systems, considered in the Markovian
framework. We propose a mixed appro\-xi\-ma\-tion scheme for the following
backward stochastic variational inequality%
\[
dY_{t}+F(t,X_{t},Y_{t},Z_{t})dt+G(t,X_{t},Y_{t})dA_{t}\in\partial\varphi
(Y_{t})dt+Z_{t}dW_{t},
\]
where $(X_{t},A_{t})_{t\in\left[  0,T\right]  }$ is the unique solution of a
reflected forward stochastic diffe\-ren\-tial equation. More precisely, we use
an Euler scheme type for the system of decoupled forward-backward variational
inequality, combined with Yosida approximation techniques.

\end{abstract}

\footnotetext[1]{{\scriptsize \textit{Acknowledgements}: The work of the first
author was supported by POSDRU/89/1.5/S/49944 project and, for the second one,
by IDEI 395/2007 project.}}

\footnotetext{\textit{{\scriptsize E-mail addresses:}}
{\scriptsize lucianmaticiuc@yahoo.com (Lucian
Maticiuc),~eduard.rotenstein@uaic.ro (Eduard Rotenstein).}}

\medskip

\medskip

\section{Introduction}

The \textit{stochastic differential equations} (SDE) with reflecting boundary
conditions, also called \textit{reflected} \textit{stochastic differential
equations }(RSDE) appears from the modeling of different kinds of constrained
phenomenon. The elliptic and parabolic PDEs with Neumann type and mixed
boundary conditions lead us to probabilistic interpretations, via the
Feynman--Ka\c{c} formula, in terms of reflected diffusion processes, which are
solutions of RSDEs. This type of equations were studied for the first time by
Skorohod in \cite{S} and after this, considered in general domains (see
\cite{LS}, \cite{M}, \cite{S1}, \cite{S2}...).

Since 1990s a lot of researchers focused their attention to numerical schemes,
methods and algorithms for the study of the behavior of the solution for
RSDEs. In the recent years, some new techniques consist in splitting-step
algorithms and mixed penalization methods. The Euler approximation was
considered for the first time by Chitashvili and Lazrieva in \cite{CL},
followed by the Euler-Peano approximation, which was introduced by Saisho in
\cite{S1}. Lepingle in \cite{L} and Slominski in \cite{S2} analyzed the
corresponding numerical schemes and their rates of convergence. In order to
approximate the solution of RSDEs, the penalization method was also very
useful (see Menaldi, \cite{M}). Approximation methods for a diffusion
reflected and stopped at the boundary appear in the literature in 1998, in the
paper of Constantini, Pacchiarotti and Sartoretto \cite{CPS}. They defined a
standard Euler projected approach to stopped reflected diffusions, approach
which yields a method with weak order of convergence ($1/2$ in particular) and
they give an easy example where this convergence rate is precise. Regarding
adaptive approximations of one-dimensional reflected Brownian motion, it can
be used a simple method of two fixed step sizes chosen according to the
distance at the boundary.

In the paper \cite{AR}, Asiminoaei and R\u{a}\c{s}canu used a mixed method
consisting in penalization and splitting-up for the study of multivalued SDE
with reflection at the boundary of the domain. The penalization method was
also used by R\u{a}\c{s}canu in \cite{R} for the study of the generalized
Skorohod problem and of its link to multivalued SDE governed by a general
maximal monotone operator (of sudifferential type). Recently, in \cite{DZ},
Ding and Zhang combined the penalization technique with the splitting-step
idea to propose some new schemes for the RSDE in the upper half space.

In 1990, Pardoux and Peng introduced in \cite{PP} the notion of nonlinear
\textit{backward stochastic differential equation} (for short, BSDE), and they
obtained the existence and uniqueness result for this kind of equation. Since
then, the interest in BSDEs has kept growing and there have been a lot of
works on that subject, both in the direction of the generalization of the
equations that appear and in constructing schemes of approximation for them.
The \textit{backward stochastic variational inequalities} (for short, BSVI)
were analyzed by Pardoux and R\u{a}\c{s}canu in \cite{PR} and \cite{PR2} (the
extension for Hilbert spaces case) by a method that consists in a penalizing
scheme, followed by its convergence.

Starting with the paper of Pardoux and Peng \cite{PP2}, have been given a
stochastic approach for the existence problem of a solution for many type of
deterministic partial differential equations (PDE for short). In \cite{PR} it
is proved, using a probabilistic interpretation, the existence for the
viscosity solution for a multivalued PDE (with subdifferential operator) of
parabolic and elliptic type. More recently, Maticiuc and R\u{a}\c{s}canu in
\cite{MR1}, prove an extended result concerning generalized type of BSDE
(including an integral with respect to an adapted continuous increasing
function and two subdifferential operators). These type of BSVI allows to
prove Feynman-Kac type formula for the representation of the solution of PVI
with mixed nonlinear multivalued Neumann-Dirichlet boundary conditions.

Even this type of the penalization approach is very useful when we deal with
multivalued backward stochastic dynamical systems governed by a
subdifferential operator, it fails for the case of a general maximal monotone
operator. This motivated a new approach, via convex analysis, for the study of
both forward and backward multivalued differential systems. In \cite{RR},
R\u{a}\c{s}canu and Rotenstein identified the solutions of those type of
equations with the minimum points of some proper, convex, lower semicontinuous
functions, defined on well chosen Banach spaces.

Euler-type approximation schemes for BSDE, and for BSDE with exit time for the
forward part of the system, were introduced by Bouchard and Touzi in \cite{BT}
and Bouchard and Menozzi in \cite{BM}. They considered the Markovian framework
of a coupled forward-backward stochastic differential system and they defined
an adapted backward Euler scheme for the strong approximation of the backward
SDE with finite stopping time horizon, namely the first exit time of the
forward SDE from a cylindrical domain. In \cite{BC}, Bouchard and Chassagneux
study the discrete-time approximation of the solution of a BSDE with a
reflecting barrier.

The paper is organized as follows. Section 2 presents some basic notations,
hypothesis and results that are used throughout this paper. Section 3 is
dedicated to the analysis of the behavior of an approximation scheme defined
for a backward stochastic variational inequality. In Section 4 we present an
existence and uniqueness result for a generalized BSVI and we propose a mixed
Euler type approximation scheme for its solution.

\section{Notations. Hypothesis. Preliminaries}

In all that follows we shall consider a finite horizon $T>0$ and a complete
probability\ space $\left(  \Omega,\mathcal{F},\mathbb{P}\right)
\mathcal{\ }$on which is defined a standard $d$-dimensional Brownian motion
$W=\left(  W_{t}\right)  _{t\leq T}$ whose natural filtration is denoted
$\mathbb{F=}\{\mathcal{F}_{t},\ 0\leq t\leq T\}.$ More precisely, $\mathbb{F}$
is the filtration generated by the process $W$ and augmented by $\mathcal{N}%
_{\mathbb{P}}$, the set of all $\mathbb{P}$-null sets, i.e. $\mathcal{F}%
_{t}=\sigma\{W_{s},$\ $s\leq t\}\vee\mathcal{N}_{\mathbb{P}}$.$\smallskip$

We denote by $L_{ad}^{r}(\Omega;C(\left[  0,T\right]  ;\mathbb{R}^{k})),$
$r\in\lbrack1,\infty),$ the closed linear subspace of adapted stochastic
processes $f\in L^{r}(\Omega,\mathcal{F},\mathbb{P};C(\left[  0,T\right]
;\mathbb{R}^{k}))$, i.e. $f(\cdot,t):\Omega\rightarrow\mathbb{R}^{k}$ is
$\mathcal{F}_{t}$-measurable for all $t\in\left[  0,T\right]  $ and
$\mathbb{E}\left(  \sup_{t\in\left[  0,T\right]  }\left\vert f\left(
t\right)  \right\vert ^{r}\right)  <\infty$. Also, we shall use the notation
$L_{ad}^{r}(\Omega;L^{q}(\left[  0,T\right]  ;\mathbb{R}^{k})),$
$r,q\in\lbrack1,\infty)$ the Banach space of $\mathcal{F}_{t}$-measurable
stochastic processes $f:\Omega\times\left[  0,T\right]  \rightarrow
\mathbb{R}^{k}$ such that $\mathbb{E}\left(  \int_{0}^{T}\left\vert f\left(
t\right)  \right\vert ^{q}dt\right)  ^{r/q}<\infty.$\medskip

\noindent Consider the following data:\vspace{-0.06in}

\begin{itemize}
\item the continuous coefficient functions $b:\left[  0,T\right]
\times\mathbb{R}^{m}\rightarrow\mathbb{R}^{m},$ $\sigma:\left[  0,T\right]
\times\mathbb{R}^{m}\rightarrow\mathbb{R}^{m\times d},$ $g:\mathbb{R}%
^{m}\rightarrow\mathbb{R}^{n}$ and $F:\left[  0,T\right]  \times\mathbb{R}%
^{m}\times\mathbb{R}^{n}\times\mathbb{R}^{n\times d}\rightarrow\mathbb{R},$
which satisfies the following standard assumptions:

for some constants $\alpha\in\mathbb{R},$ $L,$ $\beta,$ $\gamma\geq0$ and for
all $t\in\left[  0,T\right]  ,$ $x,$ $\tilde{x}\in\mathbb{R}^{m},$ $y,$
$\tilde{y}\in\mathbb{R}^{n}$ and $z,$ $\tilde{z}\in\mathbb{R}^{n\times d}:$%
\begin{equation}%
\begin{array}
[c]{cl}%
\left(  i\right)  & \left\vert b\left(  t,x\right)  -b\left(  t,\tilde
{x}\right)  \right\vert +\left\Vert \sigma\left(  t,x\right)  -\sigma\left(
t,\tilde{x}\right)  \right\Vert \leq L\left\vert x-\tilde{x}\right\vert
,\medskip\\
\left(  ii\right)  & \left\langle y-\tilde{y},F(t,x,y,z)-F(t,x,\tilde
{y},z)\right\rangle \leq\alpha\left\vert y-\tilde{y}\right\vert ^{2}%
,\medskip\\
\left(  iii\right)  & \left\vert F(t,x,y,z)-F(t,x,y,\tilde{z})\right\vert
\leq\beta\left\Vert z-\tilde{z}\right\Vert ,
\end{array}
\label{coeff assumpt}%
\end{equation}
and there exist some constants $M>0$ and $p,$ $q\in\mathbb{N}$ such that, for
all $t\in\left[  0,T\right]  $, $x\in\mathbb{R}^{m}$ and $y\in\mathbb{R}^{n}:$%
\begin{equation}
\left\vert g(x)\right\vert \leq M(1+\left\vert x\right\vert ^{q}%
)\quad\text{and}\quad\left\vert F(t,x,y,0)\right\vert \leq M(1+\left\vert
x\right\vert ^{p}+\left\vert y\right\vert ). \label{sublinear assumpt}%
\end{equation}

\item the function $\varphi:\mathbb{R}^{n}\rightarrow(-\infty,+\infty]$ which
is a proper convex lower semicontinuous function and satisfies that there
exist $M>0$ and $r\in\mathbb{N}$ such that, for all $x\in\mathbb{R}^{m}:$%
\begin{equation}
\left\vert \varphi(g(x))\right\vert \leq M(1+\left\vert x\right\vert ^{r}).
\label{fi assumpt}%
\end{equation}

\end{itemize}

The following theorem summarizes some already well known results concerning
forward and backward SDE, considered in the Markovian framework (for the proof
see Karatzas \& Shreve \cite{KS}, for forward case, and Pardoux \&
R\u{a}\c{s}canu \cite{PR} for the backward system).

\begin{theorem}
Let $\left(  t,x\right)  \in\left[  0,T\right]  \times\mathbb{R}^{m}$ be
fixed. Under the assumptions (\ref{coeff assumpt}), (\ref{sublinear assumpt})
and (\ref{fi assumpt}), the forward-backward coupled system%
\begin{equation}
\left\{
\begin{array}
[c]{l}%
dX_{s}^{t,x}=b(s,X_{s}^{t,x})ds+\sigma(s,X_{s}^{t,x})dW_{s}~,\text{ }%
s\in\left[  0,T\right]  ,\medskip\\
dY_{s}^{t,x}+F(s,X_{s}^{t,x},Y_{s}^{t,x},Z_{s}^{t,x})ds\in\partial
\varphi(Y_{s}^{t,x})ds+Z_{s}^{t,x}dW_{s}~,\text{ }t\in\left[  0,T\right]
,\medskip\\
X_{t}^{t,x}=x,\text{ }Y_{T}^{t,x}=g(X_{T}^{t,x}),
\end{array}
\right.  \label{forward_backward_system}%
\end{equation}
has a unique solution, i.e. there exist a unique process $X^{t,x}\in
L_{ad}^{2}(\Omega;C(\left[  0,T\right]  ;\mathbb{R}^{m}))$ such that%
\begin{equation}
X_{s}^{t,x}=x+\int_{t}^{t\vee s}b(r,X_{r}^{t,x})dr+\int_{t}^{t\vee s}%
\sigma(r,X_{r}^{t,x})dW_{r}~,\text{ }s\in\left[  0,T\right]  ,
\label{integral_forward_system}%
\end{equation}
and respectively%
\[
(Y^{t,x},Z^{t,x},U^{t,x})\in L_{ad}^{2}(\Omega;C(\left[  0,T\right]
;\mathbb{R}^{n}))\times L_{ad}^{2}(\Omega;L^{2}(\left[  0,T\right]
;\mathbb{R}^{n\times d}))\times L_{ad}^{2}(\Omega;L^{2}(\left[  0,T\right]
;\mathbb{R}^{n})),
\]
such that%
\begin{equation}
\left\{
\begin{array}
[c]{l}%
\displaystyle Y_{s}^{t,x}+\int_{s}^{T}U_{r}^{t,x}dr=g(X_{T}^{t,x})+\int
_{s}^{T}\mathbf{1}_{\left[  t,T\right]  }(r)F(r,X_{r}^{t,x},Y_{r}^{t,x}%
,Z_{r}^{t,x})dr-\int_{s}^{T}Z_{r}^{t,x}dW_{r}~,\;s\in\left[  0,T\right]
\medskip\\
U_{s}^{t,x}\in\partial\varphi\left(  Y_{s}^{t,x}\right)  ,\;d\mathbb{P\times
}ds\;\text{on }\Omega\times\left[  0,T\right]  .
\end{array}
\right.  \label{integral_backward_system}%
\end{equation}
Moreover, for all $p\geq2$, there exists some constant $C_{p}>0,$
$q\in\mathbb{N}^{\ast},$ such that, for all $t,$ $\tilde{t}\in\left[
0,T\right]  ,$ $x,$ $\tilde{x}\in\mathbb{R}^{n}:$%
\[%
\begin{array}
[c]{ll}%
\left(  j\right)  & \mathbb{E}\big(\sup\nolimits_{s\in\left[  0,T\right]
}\big|X_{s}^{t,x}\big|^{p}\big)\leq C_{p}(1+\left\vert x\right\vert
^{p}),\medskip\\
\left(  jj\right)  & \mathbb{E}\big(\sup\nolimits_{s\in\left[  0,T\right]
}\big|X_{s}^{t,x}-X_{s}^{\tilde{t},\tilde{x}}\big|^{p}\big)\leq C_{p}%
(1+\left\vert x\right\vert ^{p}+\left\vert \tilde{x}\right\vert ^{pq}%
)(|t-\tilde{t}|^{p/2}+|x-\tilde{x}|^{p}),\medskip\\
\left(  jjj\right)  & \mathbb{E}\big(\sup\nolimits_{s\in\left[  0,T\right]
}\big|Y_{s}^{t,x}\big|^{2}\big)\leq C(1+\left\vert x\right\vert ^{2}%
),\medskip\\
\left(  jv\right)  & \mathbb{E}\big(\sup\nolimits_{s\in\left[  0,T\right]
}\big|Y_{s}^{t,x}-Y_{s}^{\tilde{t},\tilde{x}}\big|^{2}\big)\leq C_{2}\left[
\mathbb{E}\big|g(X_{T}^{t,x})-g(X_{T}^{\tilde{t},\tilde{x}})\big|^{2}\right.
+\smallskip\\
& \quad\quad\quad\quad\left.  +\mathbb{E}%
{\displaystyle\int_{0}^{T}}
\big|\mathbf{1}_{\left[  t,T\right]  }(r)F(r,X_{r}^{t,x},Y_{r}^{t,x}%
,Z_{r}^{t,x})-\mathbf{1}_{[\tilde{t},T]}(r)F(r,X_{r}^{\tilde{t},\tilde{x}%
},Y_{r}^{t,x},Z_{r}^{t,x})\big|^{2}dr\right]  .
\end{array}
\]

\end{theorem}

\section{Approximations schemes for BSVI}

We will consider a partition of $\left[  0,T\right]  $,%
\[
\pi=\{t_{i}=ih:0\leq i\leq n\}\text{,}\;\text{with }h:=T/n,\;n\in
\mathbb{N}^{\ast},
\]
on which we approximate the solution of the backward stochastic variational
inequality (\ref{integral_backward_system}). For the numerical simulations of
the forward part, the most standard approach consists in approximating the SDE
in a proper way on each interval $\left[  t_{i},t_{i+1}\right]  $ by the
classical Euler scheme (see, e.g. Kloeden \& Platen \cite{KP}):%
\[
\left\{
\begin{array}
[c]{l}%
X_{t_{i+1}}^{h}=X_{t_{i}}^{h}+b\left(  X_{t_{i}}^{h}\right)  ~h+\sigma\left(
X_{t_{i}}^{h}\right)  \left(  W_{t_{i+1}}-W_{t_{i}}\right)  ,\;i=\overline
{0,n-1}\medskip\\
X_{0}^{h}=X_{0}.
\end{array}
\right.
\]
We remark that the above numerical scheme is easy to implement since it
requires only the simulation of $d$-independent Gaussian variables for the
Brownian increments, providing a weak error of $h$ order.

For $t\in\left[  t_{i},t_{i+1}\right]  $ let%
\[
X_{t}^{h}=X_{t_{i}}^{h}+b(X_{t_{i}}^{h})\left(  t-t_{i}\right)  +\sigma
(X_{t_{i}}^{h})\left(  W_{t}-W_{t_{i}}\right)  .
\]
We have the following estimation of the error given by the Euler scheme (see
\cite{KP}).

\begin{proposition}
\label{Euler for X}Under the assumptions (\ref{coeff assumpt}) on the
coefficients $b$ and $\sigma$, for all $p\geq1$, there exists $C_{p}>0$ such
that%
\[
\max_{\overline{0,n-1}}\mathbb{E}\bigg(\sup\limits_{t\in\left[  0,T\right]
}\big|X_{t}-X_{t}^{h}\big|^{p}+\sup\limits_{t\in\left[  t_{i},t_{i+1}\right]
}\big|X_{t}-X_{t_{i}}\big|^{p}\bigg)^{1/p}\leq C_{p}~\sqrt{h}.
\]

\end{proposition}

Here and subsequently we will consider the one-dimensional BSDE case.

Using the Yosida approximation $\nabla\varphi_{\varepsilon}$ of the
multivalued operator $\partial\varphi$, with $\varepsilon=h^{a}$ and
$a\in\left(  0,1/2\right)  $ (the way of choosing this constant will be
detailed later), we deduce that the following approximate equation%
\begin{equation}
Y_{t}^{h}+\int_{t}^{T}\nabla\varphi_{h^{a}}(Y_{r}^{h})dr=g(X_{T})+\int_{t}%
^{T}F(r,X_{r},Y_{r}^{h},Z_{r}^{h})dr-\int_{t}^{T}Z_{r}^{h}dW_{r},\;\forall
t\in\left[  0,T\right]  ,\;\mathbb{P}-a.s., \label{approximation of BSVI}%
\end{equation}
admits a unique solution $\left(  Y_{t}^{h},Z_{t}^{h}\right)  \in L_{ad}%
^{2}(\Omega;C(\left[  0,T\right]  ;\mathbb{R}))\times L_{ad}^{2}(\Omega
;L^{2}(\left[  0,T\right]  ;\mathbb{R}^{d})).$

Further, inspired by the paper of Bouchard and Touzi, \cite{BT}, let us define
an Euler type approximation for the Yosida approximation process
$Y^{\varepsilon}$. For an intuitive introduction, let $Y_{T}^{h}:=g(X_{T}%
^{\pi})$ be the initial condition, and, for $i=\overline{n-1,0}$, remark that%
\begin{equation}
Y_{t_{i}}^{h}\sim Y_{t_{i+1}}^{h}+h\left[  F(t_{i},X_{t_{i}}^{h},Y_{t_{i}}%
^{h},Z_{t_{i}}^{h})-\nabla\varphi_{h^{a}}(Y_{t_{i}}^{h})\right]  -Z_{t_{i}%
}^{h}(W_{t_{i+1}}-W_{t_{i}}); \label{intuitive consideration}%
\end{equation}
taking the conditional expectation $\mathbb{E}^{i}\left(  \cdot\right)
:=\mathbb{E}\left(  \cdot~|\mathcal{F}_{t_{i}}\right)  $, we obtain%
\[
Y_{t_{i}}^{h}\sim\mathbb{E}^{i}(Y_{t_{i+1}}^{h})+h\left[  F(t_{i},X_{t_{i}%
}^{h},Y_{t_{i}}^{h},Z_{t_{i}}^{h})-\nabla\varphi_{h^{a}}(Y_{t_{i}}%
^{h})\right]  .
\]
If we multiply (\ref{intuitive consideration}) by $W_{t_{i+1}}-W_{t_{i}}$ it
follows%
\[
Z_{t_{i}}^{h}\sim\dfrac{1}{h}\mathbb{E}^{i}(Y_{t_{i+1}}^{h}(W_{t_{i+1}%
}-W_{t_{i}})).
\]
Therefore, we propose the following implicit discretization procedure, which
define the pair $\left(  Y^{h},Z^{h}\right)  $ inductively, for $i=\overline
{n-1,0}:$%
\begin{equation}
\left\{
\begin{array}
[c]{l}%
\tilde{Y}_{T}^{h}:=g(X_{T}^{h}),\;\tilde{Z}_{T}^{h}=0,\medskip\\
\tilde{Y}_{t_{i}}^{h}:=\mathbb{E}^{i,h}(\tilde{Y}_{t_{i+1}}^{h})+h\left[
F(t_{i},X_{t_{i}}^{h},\tilde{Y}_{t_{i}}^{h},\tilde{Z}_{t_{i}}^{h}%
)-\nabla\varphi_{h^{a}}(\tilde{Y}_{t_{i}}^{h})\right]  ,\medskip\\
\tilde{Z}_{t_{i}}^{h}:=\dfrac{1}{h}\mathbb{E}^{i,h}(\tilde{Y}_{t_{i+1}}%
^{h}(W_{t_{i+1}}-W_{t_{i}})),\medskip\\
\tilde{U}_{t_{i}}^{h}:=\nabla\varphi_{h^{a}}(\mathbb{E}^{i,h}(\tilde
{Y}_{t_{i+1}}^{h})),
\end{array}
\right.  \label{approximation scheme}%
\end{equation}
where $\mathbb{E}^{i,h}\left(  \cdot\right)  :=\mathbb{E}\left(
\cdot~|\mathcal{F}_{t_{i}}^{h}\right)  $ and $\mathcal{F}_{t_{i}}^{h}%
:=\sigma(X_{t_{j}}^{h}:0\leq j\leq i).$

\begin{remark}
Observe that $\tilde{Y}_{t_{i}}^{h}$ is defined implicitly as the solution of
a fixed point problem. Since the involved functions are Lipschitz, it is well
defined. Moreover, for small values of $h>0$ it can be estimated numerically
in an accurate way.
\end{remark}

\begin{remark}
We can also use an explicit scheme to define%
\[
\tilde{Y}_{t_{i}}^{h}:=\mathbb{E}^{i,h}(\tilde{Y}_{t_{i+1}}^{h})+h\mathbb{E}%
^{i,h}\left[  F(t_{i},X_{t_{i}}^{h},\tilde{Y}_{t_{i+1}}^{h},\tilde{Z}_{t_{i}%
}^{h})-\nabla\varphi_{h^{a}}(\tilde{Y}_{t_{i+1}}^{h})\right]  .
\]
The advantage of this scheme is that it does not require a fixed point
procedure but, from a numerical point of view, adding a term in the
conditional expectation makes it more difficult to estimate. Therefore the
implicit scheme can be more tractable in practice.
\end{remark}

\begin{remark}
We have that the filtration $\mathcal{F}_{t}$ generated by the Brownian motion
coincides with the filtration generated by the diffusion process $X$, i.e.
$\mathcal{F}_{t}=\mathcal{F}_{t}^{X}$, and, from the Markov property of the
process $X^{h}$, it follows that%
\[%
\begin{array}
[c]{l}%
\mathbb{E}^{i}(\tilde{Y}_{t_{i+1}}^{h})=\mathbb{E}^{i,h}(\tilde{Y}_{t_{i+1}%
}^{h})=\mathbb{E}(\tilde{Y}_{t_{i+1}}^{h}~|X_{t_{i}}^{h}),\medskip\\
\mathbb{E}^{i}(\tilde{Y}_{t_{i+1}}^{h}(W_{t_{i+1}}-W_{t_{i}}))=\mathbb{E}%
^{i,h}(\tilde{Y}_{t_{i+1}}^{h}(W_{t_{i+1}}-W_{t_{i}}))=\mathbb{E}(\tilde
{Y}_{t_{i+1}}^{h}(W_{t_{i+1}}-W_{t_{i}})~|X_{t_{i}}^{h}).
\end{array}
\]

\end{remark}

Consider now a continuous version of (\ref{approximation scheme}). From the
martingale representation theorem there exists a square integrable process
$\tilde{Z}^{h}$ such that%
\begin{equation}
\tilde{Y}_{t_{i+1}}^{h}=\mathbb{E}^{i}(\tilde{Y}_{t_{i+1}}^{h})+\int_{t_{i}%
}^{t_{i+1}}\tilde{Z}_{s}^{h}dW_{s}, \label{martingale repres}%
\end{equation}
and, therefore, we define, for $t\in(t_{i},t_{i+1}],$%
\begin{equation}
\bar{Y}_{t}^{h}:=\tilde{Y}_{t_{i}}^{h}-\left(  t-t_{i}\right)  \left[
f(t_{i},X_{t_{i}}^{h},\tilde{Y}_{t_{i}}^{h},\tilde{Z}_{t_{i}}^{h}%
)-\nabla\varphi_{h^{a}}(\tilde{Y}_{t_{i}}^{h})\right]  +\int_{t_{i}}^{t}%
\tilde{Z}_{s}^{h}dW_{s}. \label{cont version of Y}%
\end{equation}

Obviously, we obtain that $\bar{Y}_{t_{i}}^{h}=\tilde{Y}_{t_{i}}^{h}$, and,
for the simplicity of the notation, we will continue to write $\tilde{Y}%
_{t}^{h}$ for $\bar{Y}_{t}^{h}$.

\begin{remark}
From (\ref{approximation scheme}), (\ref{martingale repres}) and the isometry
property, we notice that, for $i=\overline{0,n-1},$%
\begin{equation}
h~\tilde{Z}_{t_{i}}^{h}=\mathbb{E}^{i}(\tilde{Y}_{t_{i+1}}^{h}(W_{t_{i+1}%
}-W_{t_{i}}))=\mathbb{E}^{i}\left[  (W_{t_{i+1}}-W_{t_{i}})\int_{t_{i}%
}^{t_{i+1}}\tilde{Z}_{s}^{h}dW_{s}\right]  =\mathbb{E}^{i}\left[  \int_{t_{i}%
}^{t_{i+1}}\tilde{Z}_{s}^{h}ds\right]  . \label{connection between Zs}%
\end{equation}

\end{remark}

To approximate $Z_{t}^{h}$ we use%
\[
\bar{Z}_{t}^{h}:=\frac{1}{h}\mathbb{E}^{i}\left[  \int_{t_{i}}^{t_{i+1}}%
Z_{s}^{h}ds\right]  ,\;t\in\lbrack t_{i},t_{i+1})
\]
rather than $Z_{t_{i}}^{h}$, which is the best approximation in $L^{2}\left(
\Omega\times\left[  0,T\right]  \right)  $ of $Z^{h}$ by adapted processes
which are constant on each interval $[t_{i},t_{i+1})$ (see Lemma 3.4.2 from
Zhang \cite{Z1}):%
\[
\mathbb{E}\left[
{\displaystyle\int_{t_{i}}^{t_{i+1}}}
|Z_{s}^{h}-\bar{Z}_{t_{i}}^{h}|^{2}ds\right]  \leq\mathbb{E}\left[
{\displaystyle\int_{t_{i}}^{t_{i+1}}}
|Z_{s}^{h}-\eta|^{2}ds\right]  ,
\]
for all $\mathcal{F}_{t_{i}}$-measurable stochastic process $\eta$.

\begin{remark}
From (\ref{connection between Zs}), the definition of $\bar{Z}_{t_{i}}^{h}$
and Jensen inequality we obtain%
\begin{equation}%
\begin{array}
[c]{l}%
\displaystyle\mathbb{E}|\bar{Z}_{t_{i}}^{h}-\tilde{Z}_{t_{i}}^{h}|^{2}%
=\frac{1}{h^{2}}\mathbb{E}\left[  \mathbb{E}^{i}\int_{t_{i}}^{t_{i+1}}%
\Delta^{h}Z_{s}ds\right]  ^{2}\leq\frac{1}{h}\mathbb{E}\int_{t_{i}}^{t_{i+1}%
}\left[  \mathbb{E}^{i}\Delta^{h}Z_{s}\right]  ^{2}ds\leq\medskip\\
\displaystyle\leq\frac{1}{h}\mathbb{E}\int_{t_{i}}^{t_{i+1}}\mathbb{E}%
^{i}|\Delta^{h}Z_{s}|^{2}ds=\frac{1}{h}\int_{t_{i}}^{t_{i+1}}\mathbb{E}%
|\Delta^{h}Z_{s}|^{2}ds.
\end{array}
\label{ineq for approx Z}%
\end{equation}

\end{remark}

In order to prove an error estimate of the scheme first we use the solution
$\left(  Y_{t}^{h},Z_{t}^{h}\right)  _{t\in\left[  0,T\right]  }$ of the
approximating equation (\ref{approximation of BSVI}). The next result is a
straightforward consequence of Proposition 2.3 from Pardoux \& R\u{a}\c{s}canu
\cite{PR}.

\begin{proposition}
Under the assumptions (\ref{coeff assumpt})-(\ref{fi assumpt}), there exists
$C>0$ such that%
\begin{equation}
\sup_{t\in\left[  0,T\right]  }\mathbb{E}|Y_{t}-Y_{t}^{h}|^{2}+\mathbb{E}%
{\displaystyle\int_{0}^{T}}
|Z_{t}-Z_{t}^{h}|^{2}dt\leq C\Gamma\left(  T\right)  \hspace{1pt}h^{a},
\label{lemma 1}%
\end{equation}
where $\Gamma\left(  T\right)  :=\mathbb{E}\left[  1+\left\vert g(X_{T}%
)\right\vert ^{2}+\left\vert X_{T}\right\vert ^{r}+\int_{0}^{T}F\left(
0,X_{s}^{h},0,0\right)  ds\right]  $
\end{proposition}

We recall Theorem 3.4.3 from \cite{Z1}, applied for the solution $\left(
Y_{t}^{h},Z_{t}^{h}\right)  $ of (\ref{approximation of BSVI}). To otain a
similar conclusion we have to impose more restrictive assumptions than
(\ref{coeff assumpt}-\ref{fi assumpt}):

\begin{itemize}
\item there exists some constant $K>0$, such that\vspace{-0.1in}%
\begin{equation}%
\begin{array}
[c]{cl}%
\left(  i\right)  & \left\vert b\left(  x\right)  -b\left(  \tilde{x}\right)
\right\vert +\left\Vert \sigma\left(  x\right)  -\sigma\left(  \tilde
{x}\right)  \right\Vert \leq K\left\vert x-\tilde{x}\right\vert ,\;\forall
x,\tilde{x}\in\mathbb{R}^{m},\medskip\\
\left(  ii\right)  & |F(\xi)-F(\tilde{\xi})|\hspace{1pt}\leq K|\xi-\tilde{\xi
}|,\;\forall\xi,\tilde{\xi}\in\left[  0,T\right]  \times\mathbb{R}^{m}%
\times\mathbb{R\times R}^{d},\medskip\\
\left(  iii\right)  & \left\vert g\left(  y\right)  -g(\tilde{y})\right\vert
\leq K\left\vert y-\tilde{y}\right\vert ,\;\forall y,\tilde{y}\in\mathbb{R};
\end{array}
\label{cond coef Lipschitz}%
\end{equation}
\vspace{-0.25in}

\item the function $\varphi:\mathbb{R}^{n}\rightarrow(-\infty,+\infty]$ is a
proper convex lower semicontinuous function and there exist $M>0$ and
$r\in\mathbb{N}$ such that\vspace{-0.1in}%
\begin{equation}
\left\vert \varphi(g(x))\right\vert \leq M(1+\left\vert x\right\vert
^{r}),\;\forall x\in\mathbb{R}^{m}. \label{cond fi 2}%
\end{equation}

\end{itemize}

\begin{proposition}
\label{Lemma from Zhang}Let the assumptions (\ref{cond coef Lipschitz}) and
(\ref{cond fi 2}) be satisfied. We have the following estimate, for some
$C>0$,%
\[
\max_{i=\overline{0,n-1}}\sup_{t\in\left[  t_{i},t_{i+1}\right]  }%
\mathbb{E}|Y_{t}^{h}-Y_{t_{i+1}}^{h}|^{2}+%
{\displaystyle\sum\limits_{i=1}^{n}}
\mathbb{E}%
{\displaystyle\int_{t_{i}}^{t_{i+1}}}
|Z_{s}^{h}-\bar{Z}_{t_{i}}^{h}|^{2}ds\leq Ch,
\]
where $\bar{Z}_{t_{i}}^{h}:=\frac{1}{h}\mathbb{E}^{i}\left[  \int_{t_{i}%
}^{t_{i+1}}Z_{s}^{h}ds\right]  .$
\end{proposition}

\begin{proof}
The inequality%
\[
\max_{i=\overline{0,n-1}}\sup_{t\in\left[  t_{i},t_{i+1}\right]  }%
\mathbb{E}|Y_{t}^{h}-Y_{t_{i+1}}^{h}|^{2}\leq Ch
\]
can be obtained by classical calculus, using It\^{o}'s formula, Lipschitz
property of the coefficient functions and the bounds of the approximate
solution $\left(  Y_{t}^{h},Z_{t}^{h}\right)  _{t\in\left[  0,T\right]  }$ of
(\ref{approximation of BSVI}) (see Proposition 2.1 and 2.2 from \cite{PR}).

For the proof of the inequality%
\[%
{\displaystyle\sum\limits_{i=1}^{n}}
\mathbb{E}%
{\displaystyle\int_{t_{i}}^{t_{i+1}}}
|Z_{s}^{h}-\bar{Z}_{t_{i}}^{h}|^{2}ds\leq Ch
\]
is sufficient to recall the proof of Theorem 3.4.3 from \cite{Z1}.\hfill
\end{proof}

Using the estimates from the above Propositions we can prove the following:

\begin{proposition}
Let the assumptions (\ref{cond coef Lipschitz}) and (\ref{cond fi 2}) be
satisfied. Then there exists $C>0$ such that\vspace{-0.13in}%
\[
\sup_{t\in\left[  0,T\right]  }\mathbb{E}|Y_{t}^{h}-\tilde{Y}_{t}^{h}%
|^{2}+\mathbb{E}%
{\displaystyle\int_{0}^{T}}
|Z_{t}^{h}-\tilde{Z}_{t}^{h}|^{2}dt\leq Ch^{1-2a}.
\]

\end{proposition}

\begin{proof}
From (\ref{approximation of BSVI}) and (\ref{cont version of Y}) we deduce
that, for $i=\overline{0,n-1}$ and $t\in\left[  t_{i},t_{i+1}\right]  $,%
\begin{align*}
Y_{t}^{h}  &  =Y_{t_{i+1}}^{h}+\int_{t}^{t_{i+1}}\left[  F(s,X_{s},Y_{s}%
^{h},Z_{s}^{h})-\nabla\varphi_{h^{a}}(Y_{s}^{h})\right]  ds-\int_{t}^{t_{i+1}%
}Z_{s}^{h}dW_{s},\medskip\\
\tilde{Y}_{t}^{h}  &  =\tilde{Y}_{t_{i+1}}^{h}+\int_{t}^{t_{i+1}}\left[
F(t_{i},X_{t_{i}}^{h},\tilde{Y}_{t_{i}}^{h},\tilde{Z}_{t_{i}}^{h}%
)-\nabla\varphi_{h^{a}}(\tilde{Y}_{t_{i}}^{h})\right]  ds-\int_{t}^{t_{i+1}%
}\tilde{Z}_{s}^{h}dW_{s}.
\end{align*}
Throughout the proof let $\Delta^{h}F_{t}:=F(t,X_{t},Y_{t}^{h},Z_{t}%
^{h})-F(t_{i},X_{t_{i}}^{h},\tilde{Y}_{t_{i}}^{h},\tilde{Z}_{t_{i}}^{h})$,
$\Delta^{h}Y_{t}:=Y_{t}^{h}-\tilde{Y}_{t}^{h}$ and $\Delta^{h}Z_{t}:=Z_{t}%
^{h}-\tilde{Z}_{t}^{h}$, $t\in\left[  t_{i},t_{i+1}\right]  $.

Applying Energy equality we obtain that%
\begin{align}
\mathbb{E}|\Delta^{h}Y_{t}|^{2}+\int_{t}^{t_{i+1}}\mathbb{E}|\Delta^{h}%
Z_{s}|^{2}ds  &  =\mathbb{E}|\Delta^{h}Y_{t_{i+1}}|^{2}+2\mathbb{E}\int
_{t}^{t_{i+1}}\Delta^{h}Y_{s}~\Delta^{h}F_{s}ds\label{Energy equality}\\
&  -2\mathbb{E}\int_{t}^{t_{i+1}}\Delta^{h}Y_{s}\left(  \nabla\varphi_{h^{a}%
}(Y_{s}^{h})-\nabla\varphi_{h^{a}}(\tilde{Y}_{t_{i}}^{h})\right)  ds,\nonumber
\end{align}
We first compute $\Delta^{h}Y_{s}\cdot\left[  \Delta^{h}F_{s}-(\nabla
\varphi_{h^{a}}(Y_{s}^{h})-\nabla\varphi_{h^{a}}(\tilde{Y}_{t_{i}}%
^{h}))\right]  $ for which we use Lipschitz property of $F$ and $\nabla
\varphi_{h^{a}}:$%
\begin{equation}%
\begin{array}
[c]{l}%
\displaystyle2\mathbb{E}\int_{t}^{t_{i+1}}\Delta^{h}Y_{s}\cdot\left[
\Delta^{h}F_{s}-(\nabla\varphi_{h^{a}}(Y_{s}^{h})-\nabla\varphi_{h^{a}}%
(\tilde{Y}_{t_{i}}^{h}))\right]  ds\medskip\leq\\
\displaystyle\leq2K\mathbb{E}\int_{t}^{t_{i+1}}\left\vert \Delta^{h}%
Y_{s}\right\vert \cdot\left[  \left\vert s-t_{i}\right\vert +|X_{s}-X_{t_{i}%
}^{h}|+|Y_{s}^{h}-\tilde{Y}_{t_{i}}^{h}|+|Z_{s}^{h}-\tilde{Z}_{t_{i}}%
^{h}|+\frac{1}{h^{a}}|Y_{s}^{h}-\tilde{Y}_{t_{i}}^{h}|\right]  ds\leq
\medskip\\
\displaystyle\leq\left(  K^{2}\alpha+\beta\right)  \mathbb{E}\int_{t}%
^{t_{i+1}}|\Delta^{h}Y_{s}|^{2}ds+\frac{4}{\alpha}\mathbb{E}\int_{t}^{t_{i+1}%
}|X_{s}-X_{t_{i}}^{h}|^{2}ds+\frac{4}{\alpha}h^{3}+\medskip\\
\displaystyle+\frac{4}{\alpha}\mathbb{E}\int_{t}^{t_{i+1}}|Z_{s}^{h}-\tilde
{Z}_{t_{i}}^{h}|^{2}ds+\Big(\frac{4}{\alpha}+\frac{1}{\beta h^{2a}%
}\Big)\mathbb{E}\int_{t}^{t_{i+1}}|Y_{s}^{h}-\tilde{Y}_{t_{i}}^{h}|^{2}ds,
\end{array}
\label{Lipschitz consequence}%
\end{equation}
where $\alpha,\beta>0$ will be chosen later.

From now on, let $C>0$ be a constant independent of $h$, constant which can
take different values from one line to another.

From Proposition \ref{Euler for X} we have that there exists $C>0$ such that%
\[
\mathbb{E}|X_{s}-X_{t_{i}}^{h}|^{2}\leq2\mathbb{E}\left\vert X_{s}-X_{t_{i}%
}\right\vert ^{2}+2\mathbb{E}|X_{t_{i}}-X_{t_{i}}^{h}|^{2}\leq Ch,
\]
and, from Proposition \ref{Lemma from Zhang},%
\[
\mathbb{E}|Y_{s}^{h}-\tilde{Y}_{t_{i}}^{h}|^{2}\leq2\mathbb{E}|Y_{s}%
^{h}-Y_{t_{i}}^{h}|^{2}+2\mathbb{E}|Y_{t_{i}}^{h}-\tilde{Y}_{t_{i}}^{h}%
|^{2}\leq Ch+2\mathbb{E}|\Delta^{h}Y_{t_{i}}|^{2}.
\]

Using (\ref{ineq for approx Z})%
\[
\mathbb{E}|Z_{s}^{h}-\tilde{Z}_{t_{i}}^{h}|^{2}\leq2\mathbb{E}|Z_{s}^{h}%
-\bar{Z}_{t_{i}}^{h}|^{2}+2\mathbb{E}|\bar{Z}_{t_{i}}^{h}-\tilde{Z}_{t_{i}%
}^{h}|^{2}=2\mathbb{E}|Z_{s}^{h}-\bar{Z}_{t_{i}}^{h}|^{2}+\frac{2}{h}%
\int_{t_{i}}^{t_{i+1}}\mathbb{E}|\Delta^{h}Z_{s}|^{2}ds.
\]
Then (\ref{Lipschitz consequence}) yields%
\begin{equation}%
\begin{array}
[c]{c}%
\displaystyle A_{i}\left(  t\right)  :=\mathbb{E}|\Delta^{h}Y_{t}|^{2}%
+\int_{t}^{t_{i+1}}\mathbb{E}|\Delta^{h}Z_{s}|^{2}ds\leq\left(  K^{2}%
\alpha+\beta\right)  \mathbb{E}\int_{t}^{t_{i+1}}|\Delta^{h}Y_{s}|^{2}%
ds+B_{i},
\end{array}
\label{ineq for Gronwall}%
\end{equation}
where%
\begin{equation}%
\begin{array}
[c]{l}%
\displaystyle B_{i}:=\mathbb{E}|\Delta^{h}Y_{t_{i+1}}|^{2}+\frac{8}{\alpha
}\int_{t_{i}}^{t_{i+1}}\mathbb{E}|Z_{s}^{h}-\bar{Z}_{t_{i}}^{h}|^{2}%
ds+\frac{8}{\alpha}\int_{t_{i}}^{t_{i+1}}\mathbb{E}|\Delta^{h}Z_{s}%
|^{2}ds+\Big(\frac{4}{\alpha}+\frac{1}{\beta h^{2a}}\Big)Ch^{2}+\medskip\\
\;\;\;\;\;+2h\Big(\frac{4}{\alpha}+\frac{1}{\beta h^{2a}}\Big)\mathbb{E}%
\left\vert \Delta^{h}Y_{t_{i}}\right\vert ^{2}.
\end{array}
\label{def of B}%
\end{equation}
Using a backward Gronwall type inequality we deduce%
\[
\mathbb{E}|\Delta^{h}Y_{t}|^{2}\leq B_{i}e^{\displaystyle\int_{t}^{t_{i+1}%
}(K^{2}\alpha+\beta)ds}\leq B_{i}e^{\left(  K^{2}\alpha+\beta\right)  h}\leq
CB_{i},
\]
and, therefore,%
\[%
\begin{array}
[c]{c}%
\displaystyle A_{i}\left(  t\right)  \leq B_{i}+\left(  K^{2}\alpha
+\beta\right)  \int_{t}^{t_{i+1}}CB_{i}ds=B_{i}\left[  1+C\left(  K^{2}%
\alpha+\beta\right)  h\right]  \leq B_{i}\left[  1+Ch\right]  ,\;h\in\left(
0,1\right)  .
\end{array}
\]
The above inequality and the definition of $B_{i}$ implies%
\[%
\begin{array}
[c]{l}%
\displaystyle\left[  1-\left(  1+Ch\right)  \Big(\frac{4}{\alpha}+\frac
{1}{\beta h^{2a}}\Big)2h\right]  \mathbb{E}|\Delta^{h}Y_{t_{i}}|^{2}+\left[
1-\left(  1+Ch\right)  \frac{8}{\alpha}\right]  \int_{t_{i}}^{t_{i+1}%
}\mathbb{E}|\Delta^{h}Z_{s}|^{2}ds\leq\medskip\\
\displaystyle\leq\left(  1+Ch\right)  \left[  \mathbb{E}|\Delta^{h}Y_{t_{i+1}%
}|^{2}+\frac{8}{\alpha}\int_{t_{i}}^{t_{i+1}}\mathbb{E}|Z_{s}^{h}-\bar
{Z}_{t_{i}}^{h}|^{2}ds+Ch^{2-2a}\right]  .
\end{array}
\]
Taking $a\in\left(  0,1/2\right)  $, we can chose $h>0$ sufficiently small and
$\alpha,\beta>0$ large enough such that%
\[%
\begin{array}
[c]{c}%
C_{1}:=1-\left(  1+Ch\right)  \left(  \frac{4}{\alpha}+\frac{1}{\beta h^{2a}%
}\right)  2h>0\text{ and }C_{2}:=1-\left(  1+Ch\right)  \frac{8}{\alpha}>0
\end{array}
\]
and, therefore,%
\begin{equation}%
\begin{array}
[c]{l}%
\displaystyle C_{1}\mathbb{E}|\Delta^{h}Y_{t_{i}}|^{2}+C_{2}\int_{t_{i}%
}^{t_{i+1}}\mathbb{E}|\Delta^{h}Z_{s}|^{2}ds\leq\medskip\\
\displaystyle\leq\left(  1+Ch\right)  \left[  \mathbb{E}|\Delta^{h}Y_{t_{i+1}%
}|^{2}+\frac{8}{\alpha}\int_{t_{i}}^{t_{i+1}}\mathbb{E}|Z_{s}^{h}-\bar
{Z}_{t_{i}}^{h}|^{2}ds+Ch^{2-2a}\right]  .
\end{array}
\label{ineq for conclusion}%
\end{equation}
Writing the above inequality for each $i=\overline{0,n-1}$, we can deduce%
\[
\mathbb{E}|\Delta^{h}Y_{t_{i}}|^{2}\leq\left(  1+Ch\right)  ^{n}\left[
h^{2-2a}+\mathbb{E}|\Delta^{h}Y_{T}|^{2}+%
{\displaystyle\sum\limits_{i=1}^{n}}
{\displaystyle\int_{t_{i}}^{t_{i+1}}}
\mathbb{E}|Z_{s}^{h}-\bar{Z}_{t_{i}}^{h}|^{2}ds\right]  ,\;i=\overline
{0,n-1}.
\]
From the Lipschitz property of $g$ and Proposition \ref{Lemma from Zhang} we
obtain, for each $i=\overline{0,n-1}$, since $a\in\left(  0,1/2\right)  $,
that%
\begin{equation}
\mathbb{E}|\Delta^{h}Y_{t_{i}}|^{2}\leq Ch,\;\forall h\in\left(  0,1\right)
\;\text{small enough.} \label{estimation for delta Y}%
\end{equation}

For the proof of the inequality concerning $\left\Vert \Delta^{h}%
Z_{s}\right\Vert _{L^{2}\left(  \Omega\times\left[  0,T\right]  \right)  }$ we
act in the following manner (see, e.g. Bouchard \& Touzi \cite{BT}). From
(\ref{ineq for conclusion}) it follows that%
\[%
\begin{array}
[c]{l}%
\displaystyle\int_{0}^{T}\mathbb{E}|\Delta^{h}Z_{s}|^{2}ds=%
{\displaystyle\sum\limits_{i=0}^{n-1}}
\int_{t_{i}}^{t_{i+1}}\mathbb{E}|\Delta^{h}Z_{s}|^{2}ds\leq\medskip\\
\displaystyle\leq\left(  1+Ch\right)  \left[
{\displaystyle\sum\limits_{i=0}^{n-1}}
\mathbb{E}|\Delta^{h}Y_{t_{i+1}}|^{2}+Ch^{2-2a}~n+\frac{8}{\alpha}%
{\displaystyle\sum\limits_{i=0}^{n-1}}
\int_{t_{i}}^{t_{i+1}}\mathbb{E}|Z_{s}^{h}-\bar{Z}_{t_{i}}^{h}|^{2}ds\right]
-C_{1}%
{\displaystyle\sum\limits_{i=0}^{n-1}}
\mathbb{E}|\Delta^{h}Y_{t_{i}}|^{2}=\medskip\\
\displaystyle=\left(  1+Ch\right)  \left[  \mathbb{E}|\Delta^{h}Y_{T}%
|^{2}+Ch^{1-2a}+\frac{8}{\alpha}%
{\displaystyle\sum\limits_{i=0}^{n-1}}
\int_{t_{i}}^{t_{i+1}}\mathbb{E}|Z_{s}^{h}-\bar{Z}_{t_{i}}^{h}|^{2}ds\right]
+\medskip\\
+\left(  \left(  1+Ch\right)  +\left(  1+Ch\right)  \big(\frac{4}{\alpha
}+\frac{1}{\beta h^{2a}}\big)2h-1\right)
{\displaystyle\sum\limits_{i=1}^{n-1}}
\mathbb{E}|\Delta^{h}Y_{t_{i}}|^{2}-C_{1}\mathbb{E}|\Delta^{h}Y_{0}|^{2}.
\end{array}
\]
and, therefore, from (\ref{estimation for delta Y}),%
\[%
\begin{array}
[c]{l}%
\displaystyle\int_{0}^{T}\mathbb{E}|\Delta^{h}Z_{s}|^{2}ds\leq C\left[
\mathbb{E}|\Delta^{h}Y_{T}|^{2}+%
{\displaystyle\sum\limits_{i=0}^{n-1}}
\int_{t_{i}}^{t_{i+1}}\mathbb{E}|Z_{s}^{h}-\bar{Z}_{t_{i}}^{h}|^{2}%
ds+h^{1-2a}\right]  +\\
\displaystyle+\left(  Ch+\frac{8}{\alpha}h+\frac{2}{\beta}h^{1-2a}+\frac
{8}{\alpha}h^{2}+\frac{2}{\beta}h^{2-2a}\right)
{\displaystyle\sum\limits_{i=1}^{n-1}}
\mathbb{E}|\Delta^{h}Y_{t_{i}}|^{2}\leq\medskip\\
\displaystyle\leq C\left[  \mathbb{E}|\Delta^{h}Y_{T}|^{2}+%
{\displaystyle\sum\limits_{i=0}^{n-1}}
\int_{t_{i}}^{t_{i+1}}\mathbb{E}|Z_{s}^{h}-\bar{Z}_{t_{i}}^{h}|^{2}%
ds+h^{1-2a}+h^{1-2a}Ch~n\right]  =\medskip\\
\displaystyle=C\left[  Ch+Ch+Ch^{1-2a}\right]  \leq Ch^{1-2a}.
\end{array}
\]
Using the definition (\ref{def of B}) of $B_{i}$ we deduce that $B_{i}\leq Ch$
and, respectively, $\max\mathbb{E}|\Delta^{h}Y_{t_{i}}|^{2}\leq Ch$, which
completes the proof.\hfill\medskip
\end{proof}

Consequently we have proved our main result:

\begin{theorem}
There exists the constant $C>0$ which depends only on the Lipschitz constants
of the coefficients, such that:%
\begin{equation}
\sup_{t\in\left[  0,T\right]  }\mathbb{E}|Y_{t}-\tilde{Y}_{t}^{h}%
|^{2}+\mathbb{E}%
{\displaystyle\int_{0}^{T}}
\left[  |Y_{t}-\tilde{Y}_{t}^{h}|^{2}+|Z_{t}-\tilde{Z}_{t}^{h}|^{2}\right]
dt\leq Ch^{a\wedge\left(  1-2a\right)  }. \label{error estimate}%
\end{equation}

\end{theorem}

\section{Generalized BSVI. A proposed scheme for numerical approximation}

Let $\mathcal{D}$ be a open bounded subset of $\mathbb{R}^{d}$ of the form%

\[
\mathcal{D}=\{x\in\mathbb{R}^{d}:\ell\left(  x\right)  <0\},\;\;Bd\left(
\mathcal{D}\right)  =\{x\in\mathbb{R}^{d}:\ell\left(  x\right)
=0\},\smallskip
\]
where $\ell\in C_{b}^{3}\left(  \mathbb{R}^{d}\right)  $, $\left\vert
\nabla\ell\left(  x\right)  \right\vert =1,\;$for all $x\in Bd\left(
\mathcal{D}\right)  $. We know that (see, e.g., Lions \& Sznitman \cite{LS}),
for every $\left(  t,x\right)  \in\mathbb{R}_{+}\times\overline{\mathcal{D}}$,
there exists a unique pair of $\overline{\mathcal{D}}\times\mathbb{R}_{+}%
-$valued progressively measurable continuous processes $(X_{s}^{t,x}%
,A_{s}^{t,x})_{s\geq0},$ solution of the reflected SDE:%
\begin{equation}
\left\{
\begin{array}
[c]{l}%
\displaystyle X_{s}^{t,x}=x+\int_{t}^{s\vee t}b(r,X_{r}^{t,x})dr+\int
_{t}^{s\vee t}\sigma(r,X_{r}^{t,x})dW_{r}-\int_{t}^{s\vee t}\nabla\ell
(X_{r}^{t,x})dA_{r}^{t,x},\medskip\\
s\longmapsto A_{s}^{t,x}\text{\ is increasing,}\medskip\\
\displaystyle A_{s}^{t,x}=\int_{t}^{s\vee t}\mathbf{1}_{\{X_{r}^{t,x}\in
Bd\left(  \mathcal{D}\right)  \}}dA_{r}^{t,x}.
\end{array}
\right.  \label{defX}%
\end{equation}
Moreover, it can be proved that%
\[
\mathbb{E}\underset{s\in\left[  0,T\right]  }{\sup}\left(  \big|X_{s}%
^{t,x}-X_{s}^{t^{\prime},x^{\prime}}\big|^{p}\right)  \leq C(\big|x-x^{\prime
}\big|^{p}+\big|t-t^{\prime}\big|^{p/2})
\]
and, for all $\mu>0$, $\mathbb{E(}e^{\mu A_{T}^{t,x}})<\infty.\smallskip$

Consider now the following generalized backward stochastic variational
inequality:%
\begin{equation}
\left\{
\begin{array}
[c]{l}%
dY_{t}+F\left(  t,X_{t},Y_{t},Z_{t}\right)  dt+G\left(  t,X_{t},Y_{t}\right)
dA_{t}\in\partial\varphi\left(  Y_{t}\right)  dt+Z_{t}dW_{t},\;0\leq t\leq
T,\medskip\\
Y_{T}=g\left(  X_{T}\right)  .
\end{array}
\right.  \label{generalized BSVI}%
\end{equation}
We will suppose that the functions $F$ and $G$ satisfy the same assumption as
the generator function $F$ from the previous section. It is known (see
Maticiuc \& R\u{a}\c{s}canu \cite{MR2}) that the above equation admits a
unique solution, i.e., for all $t\in\left[  0,T\right]  $,$\;\mathbb{P}$-a.s.,%
\[
Y_{t}+\int_{t}^{T}U_{s}ds=g\left(  X_{T}\right)  +\int_{t}^{T}F\left(
s,X_{s},Y_{s},Z_{s}\right)  ds+\int_{t}^{T}G\left(  s,X_{s},Y_{s}\right)
dA_{s}-\int_{t}^{T}Z_{s}dW_{s},
\]
where $U_{t}\in\partial\varphi\left(  Y_{t}\right)  ,\;$a.e.\ on $\Omega
\times\left[  0,T\right]  .$

\begin{theorem}
Under the considered assumptions, the generalized BSVI (\ref{generalized BSVI}%
) admits a unique solution $\left(  Y_{t},Z_{t},U_{t}\right)  $ of
$\mathcal{F}_{t}$-progressively measurable processes. Moreover, for any $0\leq
s\leq t\leq T$, we have, for some positive constant $C$:%
\begin{equation}%
\begin{array}
[c]{rl}%
\left(  a\right)  & \displaystyle\mathbb{E}\left[  \int_{s}^{t}\big(\left\vert
Y_{r}\right\vert ^{2}+\left\vert \left\vert Z_{r}\right\vert \right\vert
^{2}\big)dr+\int_{s}^{t}\left\vert Y_{r}\right\vert ^{2}dA_{r}\right]
+\mathbb{E}\underset{s\leq r\leq t}{\sup}\left\vert Y_{r}\right\vert ^{2}\leq
CM_{1},\medskip\\
\left(  b\right)  & \mathbb{E}\big(\varphi\left(  Y_{t}\right)  \big)\leq
CM_{2}\quad\text{and}\quad\displaystyle\mathbb{E}\left[  \int_{s}%
^{t}\left\vert U_{r}\right\vert ^{2}dr\right]  \leq CM_{2},\medskip
\end{array}
\label{propr}%
\end{equation}
where%
\[
M_{1}=\mathbb{E}\left[  \left\vert \xi\right\vert ^{2}+\int_{0}^{T}%
\Big(\big|F\left(  s,0,0\right)  \big|^{2}ds+\big|G\left(  s,0\right)
\big|^{2}dA_{s}\Big)\right]  \quad\text{and}\quad M_{2}=\mathbb{E}%
\big(\left\vert \xi\right\vert ^{2}+\varphi\left(  \xi\right)  \big).
\]

\end{theorem}

For the generalized system considered above, we propose a mixed approximation
scheme, considering, for the simplicity of the presentation, only the case
$\varphi\equiv0$. Consider the grid of $\left[  0,T\right]  :\pi=\{t_{i}=ih,$
$i\leq n\},$ with $h:=T/n,$ $n\in\mathbb{N}^{\ast}$ and we will define
$X^{\pi}$, the approximating Euler scheme for the reflected process $X$. We
follow the paper \cite{CPS}, where the authors present the standard projected
Euler approach to stopped reflected diffusions.%
\[
\left\{
\begin{array}
[c]{l}%
X_{0}^{\pi}=x,\quad A_{0}^{\pi}=0,\smallskip\\
\hat{X}_{t_{i+1}}^{\pi}=X_{t_{i}}^{\pi}+b(t_{i},X_{t_{i}}^{\pi})(t_{i+1}%
-t_{i})+\sigma(t_{i},X_{t_{i}}^{\pi})(W_{t_{i+1}}-W_{t_{i}}),\medskip\\
\text{Taking the projection on the domain, we define}\smallskip\\
X_{t_{i+1}}^{\pi}=\bigskip\left\{
\begin{array}
[c]{ll}%
\hat{X}_{t_{i+1}}^{\pi}\text{ }, & \hat{X}_{t_{i+1}}^{\pi}\in\overline
{\mathcal{D}},\medskip\\
\Pr_{\overline{\mathcal{D}}}(\hat{X}_{t_{i+1}}^{\pi}), & \hat{X}_{t_{i+1}%
}^{\pi}\notin\overline{\mathcal{D}},
\end{array}
\right.  \quad\text{and}\\
A_{t_{i+1}}^{\pi}=\left\{
\begin{array}
[c]{ll}%
A_{t_{i}}^{\pi}\;, & \hat{X}_{t_{i+1}}^{\pi}\in\overline{\mathcal{D}}%
,\medskip\\
A_{t_{i}}^{\pi}+||\Pr_{\overline{\mathcal{D}}}(\hat{X}_{t_{i+1}}^{\pi}%
)-\hat{X}_{t_{i+1}}^{\pi}||, & \hat{X}_{t_{i+1}}^{\pi}\notin\overline
{\mathcal{D}}.\smallskip
\end{array}
\right.
\end{array}
\right.
\]

To define an approximation scheme for the generalized BSVI
(\ref{generalized BSVI}) consider $Y_{T}^{\pi}:=g(X_{T}^{\pi})$ and, for
$i=\overline{n-1,0},$ in an intuitive manner, using the notation $\Delta
A_{t_{i}}^{\pi}:=A_{t_{i+1}}^{\pi}-A_{t_{i}}^{\pi}$ and $\Delta W_{t_{i}%
}:=W_{t_{i+1}}-W_{t_{i}}~$:%
\[
Y_{t_{i}}\sim Y_{t_{i+1}}-G(X_{t_{i+1}}^{\pi},Y_{t_{i+1}})\Delta A_{t_{i}%
}^{\pi}-Z_{t_{i}}\Delta W_{t_{i}}~.
\]

\noindent We take the conditional expectation $\mathbb{E}^{\mathcal{F}_{i}}$
and it follows%
\[
Y_{t_{i}}\sim\mathbb{E}^{\mathcal{F}_{i}}(Y_{t_{i+1}})-\mathbb{E}%
^{\mathcal{F}_{i}}[G(X_{t_{i+1}}^{\pi},Y_{t_{i+1}})\Delta A_{t_{i}}^{\pi}].
\]

\noindent This suggest to define the following approximation scheme:%
\begin{equation}
\left\{
\begin{array}
[c]{l}%
Y_{t_{i}}^{\pi}:=\mathbb{E}^{\mathcal{F}_{i}}[Y_{t_{i+1}}^{\pi}-G(X_{t_{i+1}%
}^{\pi},Y_{t_{i+1}})\Delta A_{t_{i}}^{\pi}],\quad Y_{T}^{\pi}:=g(X_{T}^{\pi
}),\smallskip\\
Z_{t_{i}}^{\pi}:=\dfrac{1}{h}\mathbb{E}^{\mathcal{F}_{t_{i}}}[Y_{t_{i+1}}%
^{\pi}\Delta W_{t_{i}}-G(X_{t_{i+1}}^{\pi},Y_{t_{i+1}})\Delta A_{t_{i}}^{\pi
}\Delta W_{t_{i}}].
\end{array}
\right.  \label{proposed scheme}%
\end{equation}

\begin{problem}
The proof of the convergence for the scheme defined by (\ref{proposed scheme})
can provide a useful tool for the approximation of the solution for the
Generalized BSVI (\ref{generalized BSVI}). For the moment this is still an
open problem, which the interested reader is kindly invited to approach it.
\end{problem}

\end{document}